\documentclass[12pt,leqno]{amsart}
\usepackage{times}
\usepackage{amsmath,amssymb,amscd,amsthm}
\usepackage{latexsym}
\usepackage{epsfig}

\newtheorem{theorem}{Theorem}
\newtheorem{lemma}{Lemma}
\newtheorem{proposition}{Proposition}

\newtheorem{corollary}{Corollary}

\newtheorem{claim}{Claim}

\newcommand{\norm}[2]{{\left\| #1 \right\|}_{#2}}
\newcommand{\f}[2]{\frac{#1}{#2}}

\newcommand{\al}{\alpha}

\newcommand{\ga}{\gamma}

\newcommand{\De}{\Delta}

\newcommand{\ka}{\kappa}
\newcommand{\la}{\lambda}

\newcommand{\La}{\Lambda}

\newcommand{\vp}{\varphi}

\newcommand{\rn}{{\mathbf R}^n}

\newcommand{\rone}{\mathbf R^1}
\newcommand{\rtwo}{\mathbf R^2}

\newcommand{\cs}{\mathcal S}

\newcommand{\cz}{\mathcal Z}

\newcommand{\intl}{\int\limits}

\newcommand{\suml}{\sum\limits}
\newcommand{\supl}{\sup\limits}

\newcommand{\p}{\partial}

\newcommand{\beq}{\begin{equation}}
\newcommand{\eeq}{\end{equation}}
\newcommand{\beqna}{\begin{eqnarray*}}
\newcommand{\eeqna}{\end{eqnarray*}}
\newcommand{\beqn}{\begin{equation*}}
\newcommand{\eeqn}{\end{equation*}}
\newcommand{\bp}{\begin{proof}}
\newcommand{\ep}{\end{proof}}
\newcommand{\bprop}{\begin{proposition}}
\newcommand{\eprop}{\end{proposition}}
\newcommand{\bt}{\begin{theorem}}
\newcommand{\et}{\end{theorem}}
\newcommand{\bex}{\begin{Example}}
\newcommand{\eex}{\end{Example}}
\newcommand{\bc}{\begin{corollary}}
\newcommand{\ec}{\end{corollary}}
\newcommand{\bcl}{\begin{claim}}
\newcommand{\ecl}{\end{claim}}
\newcommand{\bl}{\begin{lemma}}
\newcommand{\el}{\end{lemma}}
\newcommand{\deal}{(-\De)^{\alpha}}

\begin{document}

\title
[Global solutions for  quasi-geostrophic equation]
{Global well-posedness  for the 2 D quasi-geostrophic  equation
 in a critical Besov space}

\author{Atanas Stefanov}

\address{
Department of Mathematics,
University of Kansas,
Lawrence, KS~66045, USA}
\date{\today}

\thanks{Supported in part by  NSF-DMS 0300511}



\begin{abstract}
We show that the the 2 D quasi-geostrophic equation has global and 
unique strong solution, when the (large) data belongs in the critical, scale 
invariant  space $\dot{B}^{2-2\al}_{2, \infty}\cap L^{2/(2\al-1)}$.  
\end{abstract}

\maketitle
\date{today}

\section{Introduction}
In this paper, we are concerned with the mathematical properties 
of the Cauchy problem for the 
quasi-geostrophic equation in two spatial dimensions
\begin{equation}
\label{eq:1}
\left|
\begin{array}{l}
\theta_t +\ka\deal \theta + J(\theta) \nabla \theta=0 \\
\theta (0,x)=\theta^0(x),
\end{array}\right.
\end{equation}
where $J(\theta)=(-R_2\theta,  R_1\theta)$ and $\al\in[0,1]$. 
The physical meaning and the derivation 
of \eqref{eq:1} has been discussed extensively in the 
literature, we refer the interested reader to the classical book of Pedlosky, 
\cite{P}. 

Depending on the value of the parameter $\al$, one distinguishes 
between the subcritical case, $\al>1/2$, the critical case $\al=1/2$ and the 
supercritical case $\al<1/2$. It is known that the critical case 
$\al=1/2$ is especially relevant from a physical point of view, as it is 
 direct analogue of the 3 D Navier-Stokes equations. On the other hand, 
considering the family of equations \eqref{eq:1} with $\al\in[0,1]$ 
allows us  to 
understand better the influence 
of the diffusion on the evolution. 

An important scale invariance, 
associated with \eqref{eq:1} is that 
$\theta^\la(t,x)=\la^{2\al-1} \theta(\la^{2\al}t, \la x)$ 
is a solution, if $\theta$ is. It follows  that the 
space $\dot{H}^{2-2\alpha}(\rtwo)$ is critical for the problem at hand. 
A heuristic argument can be made that a 
 well-posedness theory for initial data in $H^s$, $s<2-2\alpha$ 
should not hold. Thus, we concentrate our attention to the case 
$s\geq 2-2\alpha$.

 The theory for existence of solutions and their uniqueness  vary greatly, according to the criticality of the 
index $\al$. For the crtical and supercritical case, the 
question has been studied in \cite{CL}, \cite{CC}, \cite{Ju1}, \cite{Ju2}, 
\cite{Wu}, \cite{Wu1} among others. 
 The results are that when the data is large and belongs to $H^s$, 
$s>2-2\alpha$, then one has at least a local solution, which may blow up 
after finite time. For small data in the critical space 
(or some Besov variant), Chae-Lee, \cite{CL} and then  
J. Wu, \cite{Wu}, \cite{Wu1}  have been able to show existence of 
global solutions. 

We would like to mention that the majority of these results have been 
subsequently refined to include Besov spaces of initial data with the same level 
of regularity and scaling as the corresponding Sobolev spaces. Also, various 
uniqueness and blow-up criteria have been developed, 
see for example Section \ref{sec:2} below.  However, the fundamental  question for 
existence of global, smooth solutions in the supercritical case remains open. 
We note that very recently,  in the critical case $\al=1/2$, 
Kiselev, Nazarov and Volberg, \cite{KNV} have shown the 
existence of global and smooth\footnote{The smoothness assumption in \cite{KNV} is essentially 
at the level of $H^2(\rtwo)$, while the the critical case, the critical 
Sobolev space is $H^1(\rtwo)$.} solutions for any smooth (large) initial data. 

In the subcritical case, $\al>1/2$, which is of main concern for us, the quasi-geostrophic 
equation is better understood. Local and global well-posedness results, 
as well as $L^p$ decay estimates for the solution has been shown. 

To summarize the latest results, Constantin and Wu, \cite{CW1} have shown global 
well-posedness for (the inhomogeneous version) of \eqref{eq:1}, whenever the 
data is in $H^{s}: s>2-2\al$. For small data, there are plethora of results, which we will not 
review here, since we are primarily interested in the large data regime. 
On the other hand, time-decay estimates for $\norm{\theta(t)}{L^p}$ have been shown 
in \cite{CW1} and  \cite{Ju1}, see Section \ref{sec:2} below  for further details. Finally, we mention a 
local well-posedness result for large data in $H^{2-2\al}\cap L^2$, due to Ning Ju, \cite{Ju4}. Note that the space $H^{2-2\al}$ is not scale invariant (due to the $L^2$ part of it) and thus, such solutions cannot be rescaled to global ones.

In this work, we show that the quasi-geostrophic equation is globally well-posed in the critical 
 space $\dot{B}^{2-2\al}_{2, \infty}\cap L^{2/(2\al-1)}$, that is whenever the data $\theta^0$ 
belongs to the space, there is a global and unique\footnote{For the uniqueness one has to assume in addition $\theta^0\in L^2(\rtwo)$}  solution in the same space. 
\begin{theorem}
\label{theo:1}
Let $\al\in (1/2,1)$. 
Then for any 
initial data $\theta^0\in \dot{B}^{2-2\al}_{2, \infty}(\rtwo)\cap L^{2/(2\al-1)}(\rtwo)$, 
the quasi-geostrophic equation \eqref{eq:1} has a global solution 
$$
\theta\in L^\infty([0,\infty); \dot{B}^{2-2\al}_{2, \infty}(\rtwo)\cap L^{2/(2\al-1)}(\rtwo) )
$$ 
Moreover, the solution satisfies  the {\it a priori} estimate
\begin{equation}
\label{eq:3}
\norm{\theta(t)}{\dot{B}^{2-2\al}_{2, \infty}
\cap L^{2/(2\al-1)}}\leq C_{\ka, \al}(\norm{\theta^0}{\dot{B}^{2-2\al}_{2, \infty}\cap L^{2/(2\al-1)}}+ 
 \norm{\theta^0}{L^{2/(2\al-1)}}^{M(\al)}),
\end{equation}
for all $t>0$ and $M(\al)=\max(2, 1/(2\al-1))$.  In particular, the norms remain 
bounded for $0<t<\infty$. 

In addition, if $\theta_0\in L^2(\rtwo)$, then $\theta\in L^2((0,\infty), 
H^{\al}(\rtwo))$, in fact 
\begin{equation}
\label{eq:890}
\norm{\theta}{L^2, 
H^{\al}(\rtwo))}\leq \norm{\theta^0}{L^2(\rtwo)}.
\end{equation}
For a fixed $T>0$, 
the  solution is  unique 
class of weak solutions on $[0,T]$ 
 satisfying  $\theta\in L^\infty([0,T], L^2(\rtwo))\cap L^\infty([0,T], L^{2/(2\al-1)})\cap L^2((0,T), H^{\al}(\rtwo))$.  
\end{theorem}
Several remarks are in order.
\begin{enumerate}
\item Note that 
 global solutions exist and are unique in the space $\dot{B}^{2-2\al}_{2, \infty}(\rtwo)
\cap L^{2/(2\al-1)}(\rtwo)$, when the data is in the same scale invariant 
space. Note that such space properly contains $\dot{H}^{2-2\al}(\rtwo)$.  In other words, 
taking data in $\dot{H}^{2-2\al}(\rtwo)$ guarantees the existence of 
global solution, but by \eqref{eq:3} we only know that the slightly smaller norm 
$\|\theta(t)\|_{\dot{B}^{2-2\al}_{2, \infty}\cap L^{2/(2\al-1)}}$ stays bounded. 
\item It is an interesting question, whether  Theorem \ref{theo:1} and 
more precisely \eqref{eq:3} holds in  in the case of the Sobolev space $\dot{H}^{2-2\al}$ or even 
for some Besov space in the form $B^{2-2\al}_{2, r}$ for some $r<\infty$.  We note that the main difficulty  is the estimate \eqref{eq:3} for smooth solutions, since once \eqref{eq:3} is established,  one easily deduce the global existence and uniqueness by standard arguments. 
\item The results in Theorem \ref{theo:1} apply as stated for the case of the quasi-geostrophic equation 
on the real line, but one can recover the exact same statement, if one considers \eqref{eq:1} on 
${\mathbf T}^2$. We omit the details, as they  amount to 
 a minor modification of the  proof presented below.  
\end{enumerate}
{\bf Acknowledgement:} I would like to thank Ning Ju for several stimulating 
discussions on the topic. 

\section{Preliminaries}
\label{sec:2}

\subsection{The 2 D quasigeostrophic equation - 
existence and maximum principles}
We start this section by the Resnick's theorem, \cite{R} for 
existence of weak solutions. That is whenever $\theta^0\in L^2(\rtwo)$ and 
for any $T>0$, there exists a function \\
$\theta\in L^\infty([0,T], L^2(\rtwo))\cap L^2[[0,T], H^\al(\rtwo))$, so that for any test function $\vp$
$$
\int_{\rtwo} \theta(T) \vp(T) -
\int_0^T \int_{\rtwo} \theta (J(\theta)\nabla \vp)+\ka 
\int_0^T \int_{\rtwo} ((-\De)^{\al/2} \theta) ((-\De)^{\al/2} \vp) =
\int \theta^0\vp(0,x).
$$

In his dissertation, \cite{R}, Resnick also 
established {\it the maximum principle for $L^p$ norms}, 
that is for smooth solutions of \eqref{eq:1} and $1\leq  p< \infty$, one has 
\begin{equation}
\label{eq:2}
\norm{\theta(t)}{L^p(\rtwo)}\leq \norm{\theta^0}{L^p(\rtwo)}.
\end{equation}
This was later generalized by Constanin-Wu, \cite{CW1}, \cite{CW2} 
 for the case $p=2$ and by 
C\'ordoba-C\'ordoba, \cite{CC} in the case $p=2^n$ and N. Ju, 
\cite{Ju2} for all $p\geq 2$ to actually imply a power  
rate of decay for 
$\norm{\theta(t)}{L^p(\rtwo)}$ and an exponential rate of decay, 
when one considers the equation \eqref{eq:1}
 on the torus ${\mathbf T^2}$. In the sequel, we use primarily  
 \eqref{eq:2}, but is nevertheless interesting question to determine the 
optimal rates of decay for these norms. Note that Constantin and Wu 
have shown in \cite{CW1},  that the optimal 
rate for $\norm{\theta(t)}{L^2(\rtwo)}$ is $<t>^{-1/2\al}$.  Ning Ju  
has proved  in \cite{Ju2},  that\footnote{For example, $(p-2)/2p\al\to 0$ as $p\to 2$, whereas the optimal rate is $(2\al)^{-1}$, as shown by Constantin and Wu. On the other hand, we must note that the rate of $L^p$ decay obtained by Ning Ju holds under the assumption that $\theta_0\in L^2(\rtwo)$, 
while Constantin-Wu assume that $\theta_0\in L^1(\rtwo)$.}  $\norm{\theta(t)}{L^p(\rtwo)}\leq C(\norm{\theta^0}{L^p}) 
(1+t)^{-(p-2)/2p\al}$.

\subsection{The uniqueness theorem of Constantin-Wu}
Recall the uniqueness theorem of \\ Constantin-Wu 
(Theorem 2.2, in \cite{CW1})
\begin{theorem}(Constantin-Wu)
\label{theo:2}
Assume that $\al\in (1/2, 1]$ and $p,q$ satisfy $p\geq 1, 
q>1$ and $1/p+\al/q=\al-1/2$. Then for every $T>0$, there is at most one 
weak solution of \eqref{eq:1} in $[0,T]$, satisfying 
$$
\theta\in L^\infty([0,T], L^2(\rtwo))\cap L^2[[0,T], H^\al(\rtwo))\cap 
L^q([0,T], L^p(\rtwo)).
$$
\end{theorem}
In particular, one can take $q=\infty, 1/p=\al-1/2$ to 
obtain uniqueness for weak solutions satisfying 
$\theta\in L^\infty([0,T], L^p(\rtwo))$.

\subsection{Some Fourier Analysis} 
Define the Fourier transform  by 
$$
\hat{f}(\xi)=\intl_{\rn} f(x)e^{-i x\cdot  \xi} dx
$$ 
and its inverse by
$$
 f(x)=(2\pi)^{-n} \intl_{\rn} \hat{f}(\xi)e^{i x\cdot  \xi} d\xi. 
$$ 
For a
positive, smooth and  even function $\chi:\rone\to \rone$, supported in 
 $\{\xi:|\xi|\leq 2\}$ and $\chi(\xi)=1$ for all
$|\xi|\leq 1$.  
 Define $\vp(\xi)=\chi(\xi)-\chi(2\xi)$, 
which is supported in the annulus $1/2\leq |\xi|\leq 2$. Clearly 
$\sum_{k\in \cz} \vp(2^{-k} \xi)=1$ for all $\xi\neq 0$. 

The $k^{th}$ Littlewood-Paley projection is 
$\widehat{P_k f}(\xi)=\vp(2^{-k}\xi) \hat{f}(\xi)$.  Similarly 
$P_{<k}=\suml_{l\leq k}P_l$ given by the multiplier $\chi(2^{-k}\xi)$ etc. 
Note that the kernels of $P_k,P_{<k}  $ are  uniformly  integrable and thus 
$P_k, P_{<k}:L^p\to L^p$ for $1\leq p\leq \infty$ and  
$\norm{P_k}{L^p\to L^p}\leq 
C \norm{\hat{\chi}}{L^1}$. In particular, the  bounds are
independent of $k$. 

The  kernels  of $P_k$ are smooth and real-valued\footnote{Thus for 
a real valued function $\psi$, $P_k \psi$ is real-valued as well.} 
and $P_k$ commutes with differential operators.  We will frequently use the notation 
$\psi_k(x)$ instead of  $P_k \psi$, when this will not create confusion.  

It is convenient to define the (homogeneous and inhomogeneous) 
Sobolev norms in terms of 
the Littlewood-Paley operators. Namely for any   $s\geq 0$, define for every Schwartz 
function $\psi$ th norms 
\begin{eqnarray*}
& & \norm{\psi}{\dot{H}^s}:=\left(\sum_{k=-\infty}^{\infty}
 2^{2ks} \norm{\psi_k}{L^2}^2\right)^{1/2}\\
& & \norm{\psi}{H^s}:=\left(\norm{\psi}{L^2}^2+ \sum_{k=0}^{\infty}
 2^{2ks} \norm{\psi_k}{L^2}^2\right)^{1/2}
\end{eqnarray*}
and the corresponding spaces are then obtained as the closure of 
the set of all Schwartz functions in these norms. 
 Clearly $H^s=L^2\cap \dot{H}^s$. 

Introduce the operator $\La$ acting via 
$\widehat{\La \psi}(\xi):=|\xi| \hat{\psi}(\xi)$. Clearly, by the uniform boundedness 
of $P_k$ in the scale of  $L^p$ spaces, $\norm{\La^s \psi_k}{L^p}\sim 
2^{ks} \norm{\psi_k}{L^p}$.

Next, we introduce some basic facts from the 
theory of the paraproducts, which will be useful for us, 
when estimating the contribution of the nonlinearity. \\ 
Write for any two Schwartz functions $f, g$ and any integer $k$, 
\begin{eqnarray*}
& & 
P_k(f g)=P_k(\sum_{l_1, l_2} f_{l_1} g_{l_2})= 
P_k(\sum_{l_1, l_2: |l_1-l_2|\leq 3} f_{l_1} g_{l_2})+ 
P_k(\sum_{l_1, l_2: |l_1-l_2|> 3} f_{l_1} g_{l_2})
\end{eqnarray*}
But $$
P_k(\sum_{l_1, l_2: |l_1-l_2|\leq 3} f_{l_1} g_{l_2})
= P_k(\sum_{l_1, l_2: |l_1-l_2|\leq 3, \min(l_1, l_2)>k-3} f_{l_1} g_{l_2})
$$
since by the properties of the convolution 
$2^{l_1+1}+2^{l_2+1}$ must be at least $2^{k-1}$ and 
$$
P_k(\sum_{l_1, l_2: |l_1-l_2|> 3} f_{l_1} g_{l_2})= 
P_k(\sum_{l_1, l_2: |l_1-l_2|> 3, |\max(l_1, l_2)-k|\leq 3} f_{l_1} g_{l_2})
$$
since otherwise $supp \widehat{f_{l_1} g_{l_2}}\subset 
\{\xi: |\xi|\sim 2^{\max(l_1, l_2)}\}$, which would be away from the set 
$\{\xi: |\xi|\sim 2^{k}\}$ and thus $P_k(f_{l_1} g_{l_2})=0$. 

All in all, 
\begin{equation}
\label{eq:10}
\begin{array}{c}
 P_k(f g) = P_k(\suml_{l=k-3}^\infty P_{l}f  P_{l-3\leq \cdot\leq l+3} g) +\\
   P_k(\suml_{j=-3}^3 P_{k+j} f P_{<k+j-3} g)+ P_k(\suml_{j=-3}^3 P_{k+j} g P_{<k+j-3} f).
   \end{array}
\end{equation}
We will refer to the first term as ``high-high interaction'' term, while 
the second and the third terms represent  the ``high-low interaction'' term. 
We have the following lemma, which is an application of the representation formula 
\eqref{eq:10}. 
\begin{lemma}
\label{le:1}
For every $0<s\leq 1$, $2<p,q<\infty: 1/p+1/q=1/2$, there is the estimate
$$
|\int P_k \psi_k  [J(\psi)\cdot  \nabla\psi ] dx|\leq C_s 2^{k(1-s)}\norm{\psi_k}{L^p}
(\suml_{l\geq k-3} 2^{-s(l-k)} \norm{\La^s \psi_l}{L^2}) \norm{\psi}{L^q}. 
$$
for some absolute constant $C$. 
\end{lemma} 
\begin{proof}
Integration by parts and $div(J(\theta))=0$ yield 
$$
\int P_k \psi_k  [J(\psi)\cdot  \nabla\psi ] dx=-
\int \nabla \psi_k  P_k\cdot  [J(\psi) \psi ] dx
$$
At this point, by the boundedness of the Riesz transform on $L^p$, we  treat 
$J(\psi)$ as $T\psi$, where $T:L^r\to L^r$ for all $1<r<\infty$ 
 and ignore the vector structure. 
By H\"older's inequality, 
$$
|\int \nabla \psi_k  P_k [T(\psi)\cdot \psi ] dx|\lesssim 2^k \norm{\psi_k}{L^p} 
\norm{P_k [T(\psi) \psi ]}{L^{p'}}. 
$$
By  \eqref{eq:10},   
\begin{eqnarray*}
& & \norm{P_k [T(\psi) \psi ]}{L^{p'}}\leq 
\norm{\suml_{l=k-3}^\infty P_{l}T \psi  P_{l-3\leq \cdot\leq l+3} \psi}{L^{p'}} 
+ \norm{\suml_{j=-3}^3 P_{k+j} (T \psi)  P_{<k+j-3} \psi}{L^{p'}} + \\
& & + \norm{\suml_{j=-3}^3 P_{k+j} (\psi)  P_{<k+j-3} T \psi}{L^{p'}} \leq 
\suml_{l=k-3}^\infty \norm{P_l \psi}{L^2} 
\norm{P_{l-3\leq \cdot\leq l+3} \psi}{L^q}+ \\
& & +
\suml_{j=-3}^3 \norm{P_{k+j}\psi}{L^2}\norm{P_{<k+j-3} \psi}{L^q} \leq   C(\suml_{l\geq k-3} \norm{\psi_l}{L^2}) \norm{\psi}{L^q}.
\end{eqnarray*}
The Lemma follows by the observation  $\norm{\psi_l}{L^2}\sim 2^{-ls}
\norm{\La^s \psi_l}{L^2}$ and by reshufling the $2^{ks}$.

\end{proof}

\section{Proof of Theorem \ref{theo:1}}
\noindent 
The main step of the proof of Theorem \ref{theo:1} 
is the energy estimate \eqref{eq:3}. 

We start with the assumption that we are given a smooth solution 
$\theta(t,x)$, corresponding to an initial data $\theta^0$ up to time $T$ and 
we will prove \eqref{eq:3} based on it. Assume \eqref{eq:3} for a moment for such smooth solutions. 
We will show that the global existence and uniqueness follows in a standard way from 
an approximation argument and  the Constantin-Wu 
uniqueness result, Theorem \ref{theo:2}. 

Indeed, for a given initial data $\theta^0$, take an approximating sequence in 
$\dot{B}^{2-2\al}_{2, \infty}\cap L^{2/(2\al-1)}$,  
 $\{\theta^0_l\}$ of 
smooth functions (say in the Schwartz class $\cs$). By the Constantin-Wu existence result for 
data in $H^s: s>(2-2\al)$, we have  global and smooth solutions $\theta_l(t)$. In addition, 
they will satisfy the energy estimate \eqref{eq:3}.  Moreover, by the $L^p$ maximum principle, $\norm{\theta_l(t)}{L^q}\leq \norm{\theta_l(0)}{L^q}$ for all $1<q<\infty$, in particular for 
$q=2, q=2/(2\al-1)$. 

Taking weak limits will produce a weak solution 
$\theta(t)$ of \eqref{eq:1}, corresponding to initial data $\theta^0$, so that it  satisfies the energy estimate 
\eqref{eq:3} and $\norm{\theta}{L^\infty_t L^{2/(2\al-1)}}\leq \norm{\theta^0}{L^{2/(2\al-1)}}$.  This shows the existence of a weak solution with the required smoothness of the initial data.  

For the uniqueness part, we should require in addition that $\theta^0\in L^2(\rtwo)$. Then, we will show $\norm{\theta}{L^2_t H^\al_x}<\norm{\theta^0}{L^2}$, which allows us to apply  the  Constantin-Wu 
uniqueness result (Theorem \ref{theo:2}). That is, $\theta$  is the 
 unique solution in the class 
 $L^\infty([0,T], L^2(\rtwo))\cap L^2[[0,T], H^\al(\rtwo))\cap 
L^\infty([0,T], L^{2/(2\al-1)}(\rtwo))$. Thus, it remains to prove \eqref{eq:3} for smooth solutions 
and \eqref{eq:890}. Since, \eqref{eq:890} is relatively easy, we start with \eqref{eq:3}. 
\subsection{Proof of the energy estimate \eqref{eq:3}}
 Let  $s_0=2-2\al$. 
Take a Littlewood-Paley operator on both sides of \eqref{eq:1} 
\begin{eqnarray*}
\p_t \theta_k+\ka \deal \theta_k +P_k(J(\theta) \nabla \theta)=0.
\end{eqnarray*}
Taking a dot product with $\theta_k$ (which is real-valued!) yields 
\begin{eqnarray*}
\p_t \norm{\theta_k}{L^2}^2  +2 \ka \norm{(-\De)^{\al/2}\theta_k}{L^2}^2
 +2 \int P_k\theta_k J(\theta) \nabla \theta=0.
\end{eqnarray*}
By the properties of the Littlewood-Paley operators, 
$\norm{(-\De)^{\al/2}\theta_k}{L^2}^2\sim 2^{2\al k} \norm{\theta_k}{L^2}^2$. 
For the integral term, use Lemma \ref{le:1}  with $1/p=1/2-s_0/2, 1/q=s_0/2$. We have 
\begin{eqnarray*}
& & |\int P_k\theta_k J(\theta) \nabla \theta dx|\leq C 2^{k(1-s_0)} 
\norm{\theta_k(t)}{L^q}(\suml_{l\geq  k-3}2^{-s_0(l-k)} \norm{\La^s\theta_l}{L^2})
\norm{\theta(t)}{L^p}\leq \\
& &\leq C 2^{k(1-s_0)} \norm{\theta_k(t)}{L^q} \supl_{l} \norm{\La^s\theta_l}{L^2}
\norm{\theta(t)}{L^p}.
\end{eqnarray*}
By the $L^p$ maximum principle, \eqref{eq:2}, we have 
$\norm{\theta(t)}{L^p}\leq \norm{\theta^0}{L^p}$. 
Substituting everything in the equation allows us to conclude  
\begin{equation}
\label{eq:15}
\p_t \norm{\theta_k}{L^2}^2  +c \ka 2^{2k\al}\norm{\theta_k}{L^2}^2\leq 
C 2^{k(1-s_0)} \norm{\theta^0}{L^p} 
\norm{\theta_k(t)}{L^q} \supl_{l} \norm{\La^{s_0}\theta_l}{L^2}
\end{equation}
At this point,  the argument splits in two cases with a treshold value of 
$\al=3/4$. As expected, the case $3/4\leq \al<1$ proves out to be slightly 
 simpler, so we start with it. 
\subsubsection{The case $3/4\leq \al<1$}
The significance of the restriction $\al\geq 3/4$ is in the fact that 
$s_0=2-2\al\in (0,1/2]$. Therefore 
$1/q=s_0/2\leq 1/2-s_0/2=1/p$, implying $p\leq q$.  Thus, by 
 the Sobolev embedding\footnote{or more 
appropriately the Bernstein inequality}, the boundedness of $P_k$ on $L^p$ and the $L^p$ maximum principle  imply   
$$
\norm{\theta_k(t)}{L^q}\lesssim 
2^{2k(1/p-1/q)}\norm{\theta_k(t)}{L^p}\lesssim 
 2^{k(1-2s_0)}\norm{\theta_k(t)}{L^p}\lesssim 2^{k(1-2s_0)}\norm{\theta^0}{L^p}
$$
By \eqref{eq:15}, we infer 
\begin{equation}
\label{eq:316}
\p_t \norm{\theta_k}{L^2}^2  +  c \ka 2^{2k\al}\norm{\theta_k}{L^2}^2\leq 
C 2^{k(2-3 s_0)} \norm{\theta^0}{L^p}^2
 \supl_{l} \norm{\La^{s_0}\theta_l(t)}{L^2}
\end{equation}
It is a standard step now to make use of the Gronwal's inequality, namely rewrite 
\eqref{eq:316} as 
$$
\p_t (\norm{\theta_k}{L^2}^2 e^{c \ka 2^{2k\al}t})\leq C 2^{k(2-3 s_0)}
e^{ c \ka 2^{2k\al}t} \norm{\theta^0}{L^p}^2 
 \supl_{l} \norm{\La^{s_0}\theta_l(t)}{L^2}
$$
and estimate after integration 
\begin{equation}
\label{eq:17}
\norm{\theta_k(t)}{L^2}^2\leq C_\ka 2^{k(2-3 s_0-2\al)}
\norm{\theta^0}{L^p}^2 \sup_{0\leq z\leq t}  \supl_l \norm{\La^{s_0}\theta_l(z)}{L^2}+
\norm{\theta^0_k}{L^2}^2 e^{- c \ka 2^{2 k\al}t}.
\end{equation}
Note that in the formula above $C_k\sim 1/\ka$ and $2-3 s_0-2\al=-2s_0$. \\
Introduce the functional 
$$
J(t)=\sup_{0\leq z\leq t} \sup_k 2^{k s_0} \norm{\theta_k(z)}{L^2}.
$$
Clearly, one may deduce from  \eqref{eq:17} that 
$$
J^2(t)\leq J^2(0) + C_\ka  J(t) \norm{\theta^0}{L^p}^2,
$$
hence 
$$
J(t)\leq 2 J(0)+ C_\ka \norm{\theta^0}{L^p}^2,
$$
which is 
\begin{equation}
\label{eq:19}
\sup_k 2^{k(2-2\al)} \norm{\theta_k(t)}{L^2}\leq 
2 \sup_k 2^{k(2-2\al)} \norm{\theta^0_k}{L^2}+C_\ka 
\norm{\theta^0}{L^p}^2. 
\end{equation}
This is the {\it a priori} estimate of the solution $\theta$,  \eqref{eq:3} for the case $\al\in [3/4,1)$. 
As we have observed in the beginning of the section, 
 it follows that the 2 D quasi-geostrophic 
equation \eqref{eq:1} has  
global solution with (potentially large) 
data in the scale invariant space 
$\dot{B}^{2-2\al}_{2, \infty}(\rtwo)\cap L^{2/(2\al-1)}(\rtwo)$.

\subsubsection{The case $1/2< \al<3/4$}
In this case, it is clear that $s_0=2-2\al\in (1/2, 1)$, whence $2<q=(1-\al)^{-1}<p=(\al-1/2)^{-1}$. 
By Gagliardo-Nirenberg's,  
$$
\norm{\theta_k}{L^q}\leq C \norm{\La^{2-2\al}\theta_k}{L^2}^\ga
\norm{\La^{-a}\theta_k}{L^p}^{1-\ga},
$$
with $\ga=\frac{3-4\al}{2-2\al}\in (0,1)$ and 
$a=\f{(2-2\al)(3-4\al)}{2\al-1}$. Thus, by $\norm{\La^{-a}\theta_k}{L^p}\sim 
2^{-a k}\norm{\theta_k}{L^p}$, whence 
it follows 
that 
$$
\norm{\theta_k}{L^q}\leq C 2^{-k(3-4\al)}  
\sup_l \norm{\La^{s_0}\theta_l}{L^2}^\ga
\norm{\theta_k(t)}{L^p}^{1-\ga}. 
$$
Substituting this in \eqref{eq:15} yields 
\begin{equation}
\label{eq:318}
\p_t \norm{\theta_k(t)}{L^2}^2 + c  \ka 2^{2 k\al} \norm{\theta_k(t)}{L^2}^2
\leq 2^{k(1-s_0-3+4\al)}\sup_l \norm{\La^{s_0}\theta_l}{L^2}^{1+\ga}\norm{\theta_k(t)}{L^p} \norm{\theta^0}{L^p}^{1-\ga}. 
\end{equation}
Using the maximum principle $\|\theta_k(t)\|_{L^p}\lesssim \|\theta^0\|_{L^p}$, 
this reduces to 
$$
\p_t \norm{\theta_k(t)}{L^2}^2 + c  \ka 2^{2 k\al} \norm{\theta_k(t)}{L^2}^2
\leq 2^{k(1-s_0-3+4\al)}\sup_l \norm{\La^{s_0}\theta_l}{L^2}^{1+\ga}\norm{\theta^0}{L^p}^{2-\ga}. 
$$
By the Gronwall's inequality, we deduce 
\begin{equation}
\label{eq:21}
\norm{\theta_k(t)}{L^2}^2\leq 
\norm{\theta^0_k}{L^2}^2e^{-c\ka 2^{2k\al} t}+ C_\ka 2^{-2k s_0} 
\sup_{0\leq z\leq t} \sup_l \norm{\La^{s_0}\theta_l(z)}{L^2}^{1+\ga}
\norm{\theta^0}{L^p}^{2-\ga}. 
\end{equation}
By using the same energy functional $J(t)$  defined above, we conclude that 
$$
J^2(t)\leq J^2(0)+ C_k [J (t)]^{1+\ga} \norm{\theta^0}{L^p}^{2-\ga}.
$$
Since $1+\ga<2$, by Young's inequality 
$$
J^2(t)\leq J^2(0)+ \f{J^2(t)}{2}+ C_{\ka, \ga}
\norm{\theta^0}{L^p}^{(4-2\ga)/(1-\ga)}. 
$$
whence 
$$
J(t)\leq 2 J(0)+ C_{\ka, \ga}\norm{\theta^0}{L^p}^{(2-\ga)/(1-\ga)}.
$$
which is 
\begin{equation}
\label{eq:25}
\supl_k 2^{k(2-2\al)} \norm{\theta_k(t)}{L^2}\leq 
\supl_k 2^{k(2-2\al)} \norm{\theta^0_k}{L^2}+ 
C_{\ka, \ga}\norm{\theta^0}{L^p}^{(2-\ga)/(1-\ga)}.
\end{equation}
Again, this implies \eqref{eq:3} with $M(\al)=1/(2\al-1)$ 
and  the problem \eqref{eq:1} has  global solution in 
$\dot{B}^{2-2\al}_{2, \infty}(\rtwo)\cap L^{2/(2\al-1)}(\rtwo)$, when the initial data is 
taken in the same space. 

\subsection{$\theta\in L^\infty([0,\infty), L^2(\rtwo))\cap L^2((0,\infty), H^{\al}(\rtwo))$}
Both of these estimates are classical for smooth solutions, but we sketch their proofs for completeness. 

In fact,  $\theta\in L^\infty([0,\infty), L^2(\rtwo))$  follows from the maximum principle \eqref{eq:2}. 
For the second estimate, we multiply the equation by $\theta$ and integrate in $x$. We get 
\begin{eqnarray*}
& & \p_t \norm{\theta(t)}{L^2}^2+ \norm{\La^{\al} \theta(t)}{L^2}^2=-\int \theta [J(\theta)\nabla \theta] dx=0
\end{eqnarray*}
Time integartion now yields 
$$
\int_0^T \norm{\La^{\al} 
\theta(t)}{L^2}^2dt\leq \norm{\theta^0}{L^2}^2-\norm{\theta(T)}{L^2}^2 <\norm{\theta^0}{L^2}^2, 
$$
whence $\theta\in  L^2((0,\infty), H^{\al}(\rtwo))$.

\end{document}